\documentclass[12pt,notitlepage]{article}%
\usepackage{latexsym,epsfig,amsmath}
\usepackage{epsfig,graphicx}
\usepackage{amsfonts}
\usepackage{tabularx}
\usepackage{rotating}
\usepackage{harvard}
\usepackage{hyperref}
\usepackage{amsmath}
\usepackage{amssymb}
\usepackage{graphicx}%
\providecommand{\U}[1]{\protect\rule{.1in}{.1in}}
\newtheorem{theorem}{Theorem}

\newtheorem{condition}{Condition}

\newtheorem{lemma}{Lemma}

\newenvironment{proof}[1][Proof]{\textbf{#1.} }{\  \rule{0.5em}{0.5em}}

\begin{document}

\title{Joint Time Series and Cross-Section Limit Theory under Mixingale Assumptions}
\author{Jinyong Hahn\thanks{UCLA, Department of Economics, 8283 Bunche Hall, Mail
Stop: 147703, Los Angeles, CA 90095, hahn@econ.ucla.edu}\\UCLA
\and Guido Kuersteiner\thanks{University of Maryland, Department of Economics,
Tydings Hall 3145, College Park, MD, 20742, kuersteiner@econ.umd.edu}\\University of Maryland
\and Maurizio Mazzocco\thanks{UCLA, Department of Economics, 8283 Bunche Hall, Mail
Stop: 147703, Los Angeles, CA 90095, mmazzocc@econ.ucla.edu}\\UCLA}
\maketitle

\begin{abstract}
In this paper we complement joint time series and cross-section convergence
results of Hahn, Kuersteiner and Mazzocco (2016) by allowing for serial
correlation in the time series sample. The implications of our analysis are
limiting distributions that have a well known form of long run variances for
the time series limit. We obtain these results at the cost of imposing strict
stationarity for the time series model and conditional independence between
the time series and cross-section samples. Our results can be applied to
estimators that combine time series and cross-section data in the presence of
aggregate uncertainty in models with rationally forward looking
agents.\newline

Keywords:\ joint limiting distribution, data combination, stable convergence,
external validity

\end{abstract}

\newpage

\section{Introduction}

Aggregate shocks play an important role in decision making of rational agents.
The presence of aggregate shocks complicates model specification and inference
in cross-sectional settings. In Hahn, Kuersteiner and Mazzocco (2015,2019;
henceforth HKM2019) we present three economic models that illustrate the
challenges applied researchers face when trying to account for aggregate
uncertainty. Because cross-sectional samples contain no information about the
stochastic processes that generate aggregate uncertainty, we proposed to
combine cross-sectional data with long historical time series in structural
models that explicitly take into account how agents react to aggregate shocks.

Aggregate uncertainty is also considered by Rosenzweig and Udry (2019) in the
context of treatment effect estimators. In this paper we use a stylized
version of their model to illustrate how the methods proposed in HKM2019 can
be used to evaluate causal parameters under hypothetical scenarios regarding
the realization of an aggregate shock. Our method can be seen as a way to
address problems of external validity in the presence of aggregate shocks.

A rigorous foundation of our proposed empirical strategy requires central
limit theorems for the joint convergence of time series and cross-sectional
averages. The analysis of joint limits is complicated by the fact that scaling
factors resulting from linearized versions of our non-linear estimators often
depend on aggregate shocks, leading to a lack of stochastic independence
between the cross-section and time series samples. This lack of independence
does not disappear in large samples due to cross-sectional non-ergodicity with
regard to aggregate shocks, which we addressed by considering
stable\footnote{The concept of stable convergence was introduced by Renyi
(1963). It was further developed by Aldous and Eagleson (1978) and has found
numerous applications in probability theory, statistics and econometrics.}
rather than weak convergence in distribution. Our previous asymptotic results
were predicated on some martingale assumption of the time series data, which
may not be satisfied in various applications. Therefore it is useful to
develop further stable convergence results under different sets of
assumptions.\textbf{ }

This paper develops a joint time-series and cross-sectional central limit
theorem (CLT) for strictly stationary mixingale time series and, conditional
on an aggregate shock, independently distributed cross-sectional random
variables. Mixingales introduced by Gordin (1973) and McLeish (1975) are a
flexible tool to display temporal dependence in time series settings. We use
the concept to capture situations where estimators are based on moment
conditions or pseudo-likelihood functions that suffer from some form of
misspecification. The form of misspecification we have in mind is not strong
enough to affect consistency or convergence in distribution. However, because
the criterion function is not exactly a conditional mean or the score of a
correctly specified likelihood, the martingale difference properties
associated with these criteria do not apply to the estimators we have in mind
in this paper. \textbf{ }Rather, limiting variances display the type of long
run variance expressions typically found in heteroskedasticity and
autocorrelation consistent (HAC) variance estimators.

The results obtained in this paper differ both in terms of assumptions, scope
and proof strategies from our related work in Hahn, Kuersteiner and Mazzocco
(2016, HKM2016 hereafter). In that paper, we do not assume stationarity of the
time series data or independence of the cross-sectional data conditional on
aggregate shocks. On the other hand, we do impose martingale difference
conditions both in the cross-section and time series dimension. The results in
HKM2016 then are directed towards the case of correctly specified
estimators.\footnote{In addition, our (2016) paper derives functional central
limit theorems while in this paper we focus on the finite dimensional case.}

The proof strategy for finite dimensional convergence in the two papers also
differs markedly. Here, we exploit conditional independence between the time
series and cross-section and base our proofs on marginal stable convergence of
time series sample averages and almost sure conditional convergence in
distribution of cross-sectional sample averages. In HKM2016, we do not assume
conditional independence between cross-section and time series and directly
establish a joint stable martingale difference CLT by nesting the time-series
and cross-section samples in a spatial filtration similar to filtrations
proposed in Kuersteiner and Prucha (2013). The proof of the martingale
difference CLT then is based on establishing weak convergence in $L_{1}$ of
the characteristic function.\footnote{See Aldous and Eagleson (1978) for a
definition of weak convergence in $L_{1}$.} Since it is not clear how to map
cross-sectional data into a framework of strict stationarity, a direct proof
of joint convergence for the time series and cross-sections as in HKM2016 does
not seem to be obvious under the assumptions imposed in this paper.

In this paper, we get around the above-mentioned difficulty by establishing
two separate limiting results, and combining them using a conditional
independence assumption.\textbf{ }First, we establish cross-sectional stable
convergence by adapting an argument in Eagleson (1975) to show that the
conditional (on aggregate shocks) characteristic function of the
cross-sectional average converges almost surely, which implies stable
convergence. Second, we use the stable central limit theorem for strictly
stationary mixingale time series by Dedecker and Merlev\`{e}de (2002) to
establish stable convergence of our time series average.\footnote{To the best
of our knowledge Dedecker and Merlev\`{e}de (2002), and Ikeda (2017) are the
only two stable CLTs for mixingales. Ikeda's results, in particular his
assumptions about conditional variances (Ikeda (2017), Theorem 1, Assumption
(e)), are geared towards applications in finance. It is not clear that they
could easily be adapted to the applications we have in mind. We thus rely on
Dedecker and Merlev\`{e}de (2002).} Finally, we combine the two limiting
results under the assumption that the time series and cross-section samples
are independent conditional on the aggregate shock. This is done by exploiting
insights from Barndorff-Nielsen, Hansen, Lunde and Shephard (2008) who
analyzed joint stable convergence of a random sequence that converges stably
and a random sequence that converges conditionally in law.

\section{Model and Motivating Example\label{Model}}

We consider asymptotics for a combined cross-section and time series data set
where the cross-sectional units are iid and the time series units are strictly
stationary. For concreteness we assume an explicit generating mechanism that
describes the relationship between the time series variables $\nu_{t}$ and an
iid individual specific random variable $u_{i}$ such that the cross-sectional
variables $y_{i,t}$ are generated by the model\textbf{ }%
\begin{equation}
y_{i,t}=\Upsilon\left(  \nu_{t},u_{i}\right)  . \label{Cross-section_model}%
\end{equation}
We consider two samples, $\left\{  y_{i,t}:i=1,\ldots,n\right\}  $ for some
fixed $t$ and $\left\{  \nu_{t}:t=1,\ldots,\tau\right\}  $, and are interested
in the joint distribution of normalized averages of the random variables in
these two samples. Here, the function $\Upsilon$ is measurable with respect to
an underlying probability space to be defined in more detail below. The idea
behind the specification for $y_{i,t}$ is that the individual agent's behavior
can be fundamentally assumed to be a function of the two vectors $\nu_{t}$ and
$u_{i}$. We also assume strict exogeneity in the sense that the collection
$\left\{  u_{i}\right\}  _{i=1}^{\infty}$ is independent of the collection
$\left\{  \nu_{s}\right\}  _{s=1}^{\infty}$.

Then, conditional on $\nu_{t}$, the mean and variance of $y_{i,t}$ can be
written as
\begin{align*}
\mu\left(  \nu_{t}\right)   &  \equiv\int \Upsilon\left(  \nu_{t},u\right)
f_{u}\left(  u\right)  du,\\
\sigma^{2}\left(  \nu_{t}\right)   &  \equiv\int \Upsilon\left(  \nu
_{t},u\right)  ^{2}f_{u}\left(  u\right)  du-\left(  \int \Upsilon\left(
\nu_{t},u\right)  f_{u}\left(  u\right)  du\right)  ^{2},
\end{align*}
where $f_{u}\left(  \cdot\right)  $ denotes the marginal density of $u_{i}$.
In our applications, the cross sectional moments are correctly specified which
is equivalent to $\mu\left(  \nu_{t}\right)  =0$. To allow for serial
correlation, we assume that $\left\{  \nu_{s}\right\}  _{s=1}^{\infty}$ is a
mixingale. Mixingales are defined as processes whose mean conditional on past
values eventually converges to zero as the conditioning set is moved to the
more distant past. They therefore allow for serial correlation for any finite
time horizon. A more precise technical definition is given below in Section
\ref{Second}. For simplicity, we assume that the cross-sectional sample is
observed at time $t=1$ and that the time-series sample is collected starting
at $t=1.$

Below we provide an example illustrating that cross-sectional parameters of
interest may depend on parameters governing the time-series process.
Consistent estimation of these cross-sectional parameters can be achieved by
using plug-in estimates based on time-series estimates. Asymptotic
approximations to the two-step cross-sectional estimator require joint
limiting distributions of the cross-sectional and time-series estimates.
Because the cross-section depends on the time series by
(\ref{Cross-section_model}), limit results cannot be obtained for the
cross-sectional and time-series samples separately. Rather, convergence needs
to be shown to hold jointly between the cross-sectional and time-series
samples as well as an initial realization of the time-series process $\nu
_{1}.$ This is accomplished by establishing joint stable convergence of
$\left(  \frac{1}{\sqrt{n}}\sum_{i=1}^{n}y_{i,1},\frac{1}{\sqrt{\tau}}%
\sum_{s=1}^{\tau}\nu_{s}\right)  $. To prove joint stable convergence we use
an argument that combines several results. The details are laid out in the
subsequent sections.

To motivate the need for joint limiting distributions consider the following
simple stylized example of our general framework. Assume that $\tilde{y}%
_{i,t}$ is an outcome of interest, observed for individual $i$ at time time
$t.$ The individual is subject to a policy experiment $d_{i,t}.$ The outcome
$\tilde{y}_{i,1}$ in period $t=1$ is determined by the following linear
potential outcomes model
\begin{equation}
\tilde{y}_{i,1}\left(  d\right)  =\pi_{1}d+u_{i} \label{Potential_Outcomes}%
\end{equation}
where $d$ is a fixed constant and $\pi_{t}$ is the time varying causal effect
of interest. For example, when $d_{i,t}$ is binary such that $d\in\left\{
0,1\right\}  $ the model describes the potential outcomes $\tilde{y}%
_{i,1}\left(  0\right)  $ and $\tilde{y}_{i,1}\left(  1\right)  $ with causal
effect $\tilde{y}_{i,1}\left(  1\right)  -\tilde{y}_{i,1}\left(  0\right)
=\pi_{1}$. Typical treatment parameters such as the average treatment effect
or the treatment effect on the treated are defined as cross-sectional averages
$E\left[  \tilde{y}_{i,1}\left(  1\right)  -\tilde{y}_{i,1}\left(  0\right)
\right]  $ and $E\left[  \tilde{y}_{i,1}\left(  1\right)  -\tilde{y}%
_{i,1}\left(  0\right)  |d_{i,1}=1\right]  $ and coincide with $\pi_{1}.$ This
is because we abstract from treatment heterogeneity in our example such that
the parameter $\pi_{1}$ does not vary with $i$. Realized or observed values of
$\tilde{y}_{i,1}$ satisfy the constraint that $\tilde{y}_{i,1}=\tilde{y}%
_{i,1}\left(  d_{i,1}\right)  $ for realized values of the policy variable
$d_{i,1}$.

The parameter $\pi_{1}$ is identified from a random sample in the
cross-section with random assignment of individuals to treatment. Assume that
we observe a sample $\left(  \tilde{y}_{i,1},d_{i,1}\right)  $ of size
$n.$\textbf{ }Assume that (i) $\left(  d_{i,1},u_{i}\right)  $ is IID; (ii)
$d_{i,1}$ is independent of $u_{i}$; and (iii) $\Pr\left(  d_{i,1}=1\right)
=\frac{1}{2}$, which is known to the econometrician.\footnote{Knowledge of the
propensity score is used to simplify the calculation of standard errors in
Section \ref{Section_SE}.} A consistent estimator for $\pi_{1}$ is obtained
from $\hat{\pi}_{1}=\sum_{i=1}^{n}w_{i}\tilde{y}_{i,1}$ where $w_{i}%
=d_{i,1}/\sum_{i=1}^{n}d_{i,1}-\left(  1-d_{i,1}\right)  /\sum_{i=1}%
^{n}\left(  1-d_{i,1}\right)  $.

The parameter $\pi_{t}$ possibly varies with time. For example, a job training
program may be less effective in reducing unemployment during a recession than
during an economic boom. If an investigator had cross-sectional data at a
particular point in time and the goal were to make general predictions about
the efficacy of such a job training program then the evidence obtained from a
randomized study may not be informative about the promise of the program at a
given future point in time. We adopt the approach of HKM2019, and address this
problem by parametrizing the functional form of $\pi_{t}$. In HKM2019, we
motivate such parametrizations with structural economic models that explain
individual decision making by rational economic agents. Similarly, Rosenzweig
and Udry (2019) consider models of investment behavior where the return to
investment depends on the realization of an aggregate shock.

To keep the exposition simple we postulate a simple stylized parametric model
for $\pi_{t}$ in this paper. Assume that there is an aggregate economic shock
$\nu_{t}$ that affects the causal effect of the policy in a multiplicative way
such that $\pi_{t}=\left.  \beta\right/  \left(  1+\nu_{t}^{2}\right)  $.
Here, $\beta$ is interpreted as a deep structural parameter that can be used
to impute the causal effect for values of the economic shock $\nu_{t}$ that
may not have realized during the time period when the cross-sectional sample
was observed. Further assume that we have time series data $z_{s}$ that
satisfies the dynamic equation $z_{s}=\phi+\nu_{s}$ where $\phi$ is a fixed
time invariant constant and $\nu_{s}$ is an unobserved stationary mixingale.
If the parameter $\phi$ were known, the structural parameter $\beta$ could be
identified from the relationship $\beta=\pi_{1}\left(  1+\left(  z_{1}%
-\phi\right)  ^{2}\right)  $. The solution we propose in HKM2019 consists in
using a time series sample $z_{s}$ for $s=1,...,\tau$ to estimate the
parameter $\phi.$ In our simple setting this amounts to forming the time
series average $\hat{\phi}=\tau^{-1}\sum_{s=1}^{\tau}z_{s}=\phi+\tau^{-1}%
\sum_{s=1}^{\tau}\nu_{s}.$

The specification of $\beta$ captures the idea that the causal effect of
$d_{i,t}$ on $y_{i,t}$ is dampened in periods where $z_{t}$ deviates from its
long term stationary mean $\phi.$ The parameter $\phi$ could measure an
optimal growth path with deviations from it reducing the efficacy of policy
interventions or investments. Since deviations may be persistent, it is
plausible that $\nu_{s}$ exhibits serial correlation. For example, if $\nu
_{s}=\theta\nu_{s-1}+\varepsilon_{s}$ with $\varepsilon_{s}\sim N\left(
0,1\right)  $ and $\left\vert \theta\right\vert <1$ it is easy to see that
$E\left[  \nu_{s}|\nu_{s-k}\right]  =\theta^{k}\nu_{s-k}$ and thus that
$\nu_{s}$ is a mixingale according to the definition given in Section
\ref{Second}. While the parameter $\phi$ can be consistently estimated as a
simple sample average, the sample average is not the maximum likelihood
estimator in this scenario of a Gaussian autoregressive process. As a result,
the influence function of $\hat{\phi}$ is not a martingale difference
sequence, as would be the case for the maximum likelihood estimator. Thus, the
results in HKM2016, which require martingale difference sequences, cannot be
applied to the case where the investigator is using an inefficient estimator
for the parameter $\phi.$ In this paper, we develop the necessary asymptotic
theory to handle cases where the influence function of $\phi$ is a more
general mixingale process.

The estimators $\hat{\pi}_{1}$ and $\hat{\phi}$ obtained from the
cross-section and time series samples respectively can be combined to form the
estimator $\hat{\beta}=\hat{\pi}_{1}\left(  1+\left(  z_{1}-\hat{\phi}\right)
^{2}\right)  .$ In HKM2019 we are mostly interested in inference for the
parameter $\beta$. Using a Taylor series expansion one obtains%
\begin{equation}
\sqrt{n}\left(  \hat{\beta}-\beta\right)  =\left(  1+\nu_{1}^{2}\right)
\sqrt{n}\left(  \hat{\pi}_{1}-\pi_{1}\right)  -2\pi_{1}\nu_{1}\sqrt{\frac
{n}{\tau}}\sqrt{\tau}\left(  \hat{\phi}-\phi\right)  +o_{p}\left(  1\right)  .
\label{beta_expand}%
\end{equation}
Under an asymptotic sequence where $\sqrt{\frac{n}{\tau}}\rightarrow
\sqrt{\kappa}$ for $0<\kappa<\infty$ the limiting distribution of $\hat{\beta
}$ is determined by the joint limiting distribution of $\hat{\pi}_{1}$ and
$\hat{\phi}$ conditional on $\nu_{1}.$ Because both the cross-sectional and
the time series data depend on $\nu_{1},$ establishing joint convergence
requires special care. In HKM2016 we use the concept of stable convergence. In
this paper, a somewhat different approach, based on almost sure conditional
convergence is used. We also show that the latter implies the former.

In addition to the structural parameter $\beta$ we may be interested in causal
effects outside of the observed cross-sectional sample. Using $\hat{\beta}$ we
are in a position to impute the causal effects for arbitrary values of the
aggregate shock $\nu_{s}$. For a fixed value of $\nu,$ and thus indirectly,
$z=\phi+\nu,$ define the counterfactual causal effect $\pi\left(  \nu\right)
=\left.  \beta\right/  \left(  1+\nu^{2}\right)  $ which can be estimated by
\[
\hat{\pi}\left(  \nu\right)  =\left.  \hat{\beta}\right/  \left(  1+\nu
^{2}\right)  .
\]
Inference for $\hat{\pi}\left(  \nu\right)  $ is based on a Taylor expansion
and the result in (\ref{beta_expand})
\begin{equation}
\sqrt{n}\left(  \hat{\pi}\left(  \nu\right)  -\pi\left(  \nu\right)  \right)
=\frac{1+\nu_{1}^{2}}{1+\nu^{2}}\sqrt{n}\left(  \hat{\pi}_{1}-\pi_{1}\right)
-\frac{2\pi_{1}\nu_{1}}{1+\nu^{2}}\sqrt{\frac{n}{\tau}}\sqrt{\tau}\left(
\hat{\phi}-\phi\right)  +o_{p}\left(  1\right)  \label{pi_expand}%
\end{equation}
and relies on the same joint asymptotic limits as inference for $\hat{\beta}$.

\section{Limit Theory}

The analysis of estimators such as $\hat{\beta}$ or $\hat{\pi}\left(
\nu\right)  $ depends on the joint convergence of the tuple $\left(
U_{n,\tau},\zeta\right)  $ as $n,\tau\rightarrow\infty$, where $\zeta$ is any
measurable function of $\nu_{1}$, and where $U_{n,\tau}$ is defined as
$U_{n,\tau}\equiv\left(  \frac{1}{\sqrt{n}}\sum_{i=n}^{n}y_{i,1},\frac
{1}{\sqrt{\tau}}\sum_{s=1}^{\tau}\nu_{s}\right)  $. For ease of exposition we
assume that $y_{i,1}$ and $\nu_{s}$ take values in $\mathbb{R}$. Extensions to
the multivariate case using the Cramer-Wold theorem are straight forward.

This form of convergence where the component $\zeta$ does not depend on $n$ or
$\tau$ was first discussed by Renyi (1963). To formally define it, let
$\left(  \Omega,\mathcal{F},P\right)  $ be a probability space with a
sub-sigma field $\mathcal{C}$ and assume that $B$ is any $\mathcal{C}%
$-measurable event. Following Aldous and Eagleson (1978), the sequence
$U_{n,\tau}$ is said to converge $\mathcal{C}$-stably to a random variable $U$
defined on $\left(  \Omega,\mathcal{F},P\right)  $ if $\lim_{n,\tau}P\left(
U_{n,\tau}\leq x,B\right)  $ exists for a countable dense set of points $x$
and every set $B$ that is measurable with respect to $\mathcal{C}$. As pointed
out by Aldous and Eagleson (1978), stable convergence implies that for every
$B$ with $P\left(  B\right)  >0$ the conditional (on $B$) distribution of
$U_{n,\tau}$ converges to a well defined limit. Aldous and Eagleson (1978)
also show that stable convergence is equivalent to weak $L_{1}$ convergence of
the characteristic function of $U_{n,\tau}.$ Formally, let $\bar{\zeta}$ be
any bounded $\mathcal{C}$-measurable random variable and $\ell$ any fixed
vector conforming with $U_{n,\tau}.$ Then, $U_{n,\tau}$ converges weakly in
$L_{1}$ to $U$ if
\begin{equation}
E\left[  \bar{\zeta}\exp\left(  i\ell^{\prime}U_{n,\tau}\right)  \right]
\rightarrow E\left[  \bar{\zeta}\exp\left(  i\ell^{\prime}U\right)  \right]
\label{L1-Convergence}%
\end{equation}
for $i=\sqrt{-1}.$ The notion of weak $L_{1}$ convergence of the
characteristic function plays a key role in the proof of stable central limit
theorems in Hall and Heyde (1980) and Kuersteiner and Prucha (2013). The proof
strategy in these papers can be traced back to McLeish (1974). Eagleson (1975)
on the other hand proves a stable CLT for martingale differences by showing
that the conditional characteristic function converges almost surely.
Formally, this is stated as
\[
\lim_{n,\tau}E\left[  \exp\left(  i\ell^{\prime}U_{n,\tau}\right)
|\mathcal{C}\right]  =E\left[  \exp\left(  i\ell^{\prime}U\right)
|\mathcal{C}\right]  \text{ a.s.}%
\]
Using iterated expectations, which can be defined because $U_{n,\tau},$ $U$
and $\bar{\zeta}$ all live on the same probability space, it follows that
\begin{align*}
\left\vert E\left[  \bar{\zeta}\left(  \exp\left(  i\ell^{\prime}U_{n,\tau
}\right)  -\exp\left(  i\ell^{\prime}U\right)  \right)  \right]  \right\vert
&  =\left\vert E\left[  \bar{\zeta}\left(  E\left[  \exp\left(  i\ell^{\prime
}U_{n,\tau}\right)  -\exp\left(  i\ell^{\prime}U\right)  |\mathcal{C}\right]
\right)  \right]  \right\vert \\
&  \leq E\left[  \left\vert \bar{\zeta}\right\vert \left\vert E\left[
\exp\left(  i\ell^{\prime}U_{n,\tau}\right)  -\exp\left(  i\ell^{\prime
}U\right)  |\mathcal{C}\right]  \right\vert \right]  .
\end{align*}
By the fact that the characteristic function is uniformly bounded and using
the Lesbesgue convergence theorem it follows that almost sure convergence of
the conditional characteristic function implies stable convergence.

In HKM2016 we directly establish (\ref{L1-Convergence}) for a martingale
difference array by nesting $U_{n,\tau}$ into an increasing spatial
filtration. The construction used in that paper does not require stationarity
but exploits existing central limit theorems that are specific to martingale
difference sequences. In this paper, we relax the martingale difference
assumption at the cost of imposing stationarity. Another feature of our
setting in this paper not required in HKM2016 is crucial for the argument: We
assume that conditional on $\mathcal{C}$ the components $Y_{n}=\frac{1}%
{\sqrt{n}}\sum_{i=1}^{n}y_{i,1}$ and $Z_{\tau}=\frac{1}{\sqrt{\tau}}\sum
_{s=1}^{\tau}\nu_{s}$ are independent for all $n$ and $\tau.$ This allows to
multiplicatively separate $E\left[  \exp\left(  i\ell^{\prime}U_{n,\tau
}\right)  |\mathcal{C}\right]  $ into two components. We then show that the
conditional characteristic function of $Y_{n}$ converges almost surely. This
almost sure conditional convergence turns out to be the critical component in
showing that (\ref{L1-Convergence}) holds jointly. This argument is made
rigorous in the subsequent discussion.

\subsection{Conditional Cross-Sectional CLT\label{First}}

This section establishes a conditional central limit theorem for%
\[
\frac{1}{\sqrt{n}}\sum_{i=1}^{n}y_{i,1}%
\]
conditional on $\nu_{1}$. Recall that $\mu\left(  \nu_{1}\right)  $ and
$\sigma^{2}\left(  \nu_{1}\right)  $ denote the conditional mean and variance
of $y_{i,1}=\Upsilon\left(  \nu_{1},u_{i}\right)  $ given $\nu_{1}$, and that
we assume $\mu\left(  \nu_{1}\right)  =0$ in our application.

The proof follows the arguments in Eagleson (1975) and van der Vaart and
Wellner (1996, Lemma 2.9.5) by specializing the regularity conditions in
Eagleson (1975) to the conditionally independent case. Some additional
notation is required to clearly define the probability space involved in the
construction. Consider the product space $\mathbb{R}^{\infty}\times
\mathbb{R}^{\infty}$ with Borel fields\footnote{The sigma algebra
$\mathcal{B}^{\infty}$ is the smallest sigma algebra that contains all the
sets $\left\{  x\in\mathbb{R}^{\infty}|x=\left(  x_{1},x_{2}...,\right)
,x_{1}\in I_{1},x_{2}\in I_{2},...x_{m}\in I_{m}\right\}  $ where $I_{j}$ is
an interval $(a,b]$ on $\mathbb{R}$ - see Shiryaev (1995, p.146).}
$\mathcal{B}^{\infty}\times\mathcal{B}^{\infty}$ and product measure
$P=P_{\nu}\times P_{u}$. The infinite dimensional vectors $\nu=\left(  \nu
_{1},\nu_{2},...\right)  ^{\prime}$ and $u=\left(  u_{1},u_{2},...\right)
^{\prime}$ take values in $\mathbb{R}^{\infty}\times\mathbb{R}^{\infty}$.
Assume that $\left(  \nu,u\right)  $ is $\mathcal{X}$-measurable where
$\mathcal{X=B}^{\infty}\times\mathcal{B}^{\infty}$. A probability space of
this form can be constructed using Kolmogorov's existence theorem (see
Billingsley 1995, p.486). Suppose that $\nu_{1}$ is a $\mathcal{A}$-measurable
random variable, where $\mathcal{A}$ is a sub-sigma field of $\mathcal{X}$
such that $\mathcal{A}\subset\mathcal{X}$. Consider the probability space
$\left(  \mathbb{R}^{\infty}\times\mathbb{R}^{\infty},\mathcal{X},P\right)  $.
By Breiman (1992, Theorem 4.34 and A.46), as long as the sample space is a
complete separable metric space, a regular conditional distribution on
$\mathcal{X}$ given $\mathcal{A}\subset\mathcal{X}$ exists. As in Eagleson,
let $\omega^{\prime}\in\mathbb{R}^{\infty}\times\mathbb{R}^{\infty}$ and
consider the regular conditional probability denoted by $Q_{\omega^{\prime}%
}\left(  B,\mathcal{A}\right)  =Q_{\omega^{\prime}}\left(  B\right)  $. It
follows that for fixed $B\in\mathcal{X}$, $Q_{\omega^{\prime}}\left(
B,\mathcal{A}\right)  $ is a version of $P\left(  \left.  B\right\vert
\mathcal{A}\right)  $ and for fixed $\omega^{\prime}\in\mathbb{R}^{\infty
}\times\mathbb{R}^{\infty}$, $Q_{\omega^{\prime}}\left(  \cdot\right)  $ is a
probability measure on $\mathcal{X}$.

Consider the measure space $\left(  \mathbb{R}^{\infty}\times\mathbb{R}%
^{\infty},\mathcal{X},Q_{\omega^{\prime}}\right)  $ with expectation
$E_{\omega^{\prime}}$. By a Lemma in Eagleson (1975, p.558) the following
holds: Let $\mathcal{G}$ be a sub-sigma field of $\mathcal{X}$ such that
$\mathcal{A}\subseteq\mathcal{G}$.\footnote{It is easy to see that the proof
of Lemma 1 in Eagleson (1975) goes through when $\mathcal{A=G}$.} Then, for
$P$ almost all $\omega^{\prime}\in\mathbb{R}^{\infty}\times\mathbb{R}^{\infty
}$ and a random variable $Y$ with $E\left\vert Y\right\vert <\infty$ it
follows that $E_{\omega^{\prime}}\left[  \left.  Y\right\vert \mathcal{G}%
\right]  \left(  \omega\right)  =E\left[  \left.  Y\right\vert \mathcal{G}%
\right]  \left(  \omega\right)  $ $Q_{\omega^{\prime}}$-a.s.

Define $S_{n}\equiv\sum_{i=1}^{n}y_{i,1}$, and $Y_{n}\equiv\left.
S_{n}\right/  \sqrt{n}$. We now have the following result which is established
by combining arguments in the proofs of Eagleson (1975, Theorem 2) and van der
Vaart and Wellner (1996, Lemma 2.9.5).\footnote{Lemma 2.9.5 of van der Vaart
and Wellner (1996) does not directly apply to our context, because of the
conditional nature of our problem. However, our proof uses some of the
arguments in their proof.}

\begin{theorem}
\label{Cross-Section-CLT}Assume that for some $\delta>0$ it follows that
$E\left[  \left.  \left\vert y_{i,1}\right\vert ^{2+\delta}\right\vert
\mathcal{A}\right]  \leq K<\infty$ for some constant $K.$ Then, conditional on
$\mathcal{A}$, $S_{n}$ converges in law as $n\rightarrow\infty$ to a normal
distribution. In addition
\[
\lim_{n\rightarrow\infty}E\left[  e^{itY_{n}}|\mathcal{A}\right]
=e^{-\frac{1}{2}t^{2}\sigma^{2}\left(  \nu_{1}\right)  }\text{ a.s.}%
\]

\end{theorem}

\begin{proof}
See Appendix.
\end{proof}

\subsection{Stable Time Series CLT\label{Second}}

We follow Dedecker and Merlev\`{e}de (2002) in defining a stationary sequence.
Consider the probability space $\left(  \mathbb{R}^{\infty},\mathcal{B}%
^{\infty},P_{\nu}\right)  $ which is the second coordinate of $\left(
\mathbb{R}^{\infty}\times\mathbb{R}^{\infty},\mathcal{X},P\right)  .$ Let
$T:\mathbb{R}^{\infty}\rightarrow\mathbb{R}^{\infty}$ be a bijective
bimeasurable transformation preserving $P_{\nu}.$ Let $\mathcal{M}_{0}$ be a
sigma-algebra of $\mathcal{B}^{\infty}$\textbf{ }and let $\nu_{0}$ be
$\mathcal{M}_{0}$-measurable. Assume that $T$ satisfies $T^{-1}\mathcal{M}%
_{0}\supset\mathcal{M}_{0}$.\footnote{A mapping $T$ with this property is
called compressible, see Halmos (1956, p.11).} Define the non-decreasing
filtration $\mathcal{M}_{i}=T^{-i}\left(  \mathcal{M}_{0}\right)  $. Define
the sequence $\nu_{s}=\nu_{0}\left(  T^{s}\omega\right)  $ for $\omega
\in\mathbb{R}^{\infty}$. Let $I$ be an invariant set $T^{-1}I=I$. Let
$\mathcal{I}$ be the sigma-algebra of all invariant sets.

Now impose the following assumptions that correspond to Dedecker and
Merlev\`{e}de (2002, Theorem 1, s2.). Let $\left\Vert X\right\Vert _{p}%
\equiv\left(  E\left[  \left\vert X\right\vert ^{p}\right]  \right)  ^{1/p}$
for any random variable $X\in\mathbb{R}$ defined on the probability space
$\left(  \mathbb{R}^{\infty}\times\mathbb{R}^{\infty},\mathcal{X},P\right)  $.

\begin{condition}
\label{DM_s2}Let $S_{\nu,\tau}\equiv\sum_{s=1}^{\tau}\nu_{s}$ and $Z_{\tau
}\equiv\left.  S_{\nu,\tau}\right/  \sqrt{\tau}$.\newline i) The sequence
$Z_{\tau}^{2}$ is uniformly integrable.\newline ii) The sequence $\left\Vert
E\left[  Z_{\tau}|\mathcal{M}_{0}\right]  \right\Vert _{1}\rightarrow0$ as
$\tau\rightarrow\infty.$\newline iii) There exists a nonnegative
$\mathcal{M}_{0}$ measurable variable $\eta$ with $\eta\left(  \omega\right)
=\eta\left(  T\omega\right)  $ for all $\omega\in\mathbb{R}^{\infty}$ such
that $\left\Vert E\left[  \left.  Z_{\tau}^{2}-\eta\right\vert \mathcal{M}%
_{0}\right]  \right\Vert _{1}\rightarrow0$ as $\tau\rightarrow\infty.$
\end{condition}

Let $\varphi$ be a continuous function $\varphi:\mathbb{R\rightarrow
}\mathbb{R}$ such that $\left\vert \left(  1+x^{2}\right)  ^{-1}\varphi\left(
x\right)  \right\vert $ is bounded. Let $g\left(  x\right)  $ be the standard
Gaussian density. If Condition \ref{DM_s2} holds then it follows from Dedecker
and Merlev\`{e}de (2002, Theorem 1) that
\begin{equation}
\lim_{\tau\rightarrow\infty}\left\Vert E\left[  \left.  \varphi\left(
Z_{\tau}\right)  -\int\varphi\left(  x\sqrt{\eta}\right)  g\left(  x\right)
dx\right\vert \mathcal{M}_{k}\right]  \right\Vert _{1}=0\label{Timeseries_CLT}%
\end{equation}
for every positive integer $k$. In particular, it holds for $\mathcal{M}_{1}$,
which can be taken to be equivalent to $\mathcal{A}$ in the previous section.
The next Lemma formally establishes a relationship between
(\ref{Timeseries_CLT}) and the $L_{1}$ based definition of stable convergence
given in Aldous and Eagleson (1978).

\begin{lemma}
\label{Lemma_DM_stable}Assume that Condition \ref{DM_s2} holds. Then,
$Z_{\tau}\rightarrow_{d}\sqrt{\eta}\xi_{\nu}$ ($\mathcal{A}$-stably) where
$\eta$ is $\mathcal{A}$-measurable and $\xi_{\nu}$ is standard Gaussian
independent of $\mathcal{A}$.
\end{lemma}

\begin{proof}
See Appendix.
\end{proof}

Low level conditions for Condition \ref{DM_s2} can be given using the results
of Dedecker and Doukhan (2003). Let $\left\{  \nu_{s}\right\}  _{s\geq0}$ be a
sequence of real valued stationary random variables. Define the mixingale
coefficient $\gamma_{s}\equiv\left\Vert E\left[  \nu_{s}|\mathcal{M}%
_{0}\right]  -E\left[  \nu_{s}\right]  \right\Vert _{1}$. The definition of
$\gamma_{t}$ goes back to Gordin (1973) and for an analogous definition based
on the $L_{2}$ norm to McLeish (1975). Modified to account for stationarity
and represented in terms of $L_{1}$ norms, McLeish's definition of a mixingale
can be stated as follows: A sequence $\left(  \nu_{s},\mathcal{M}_{s}\right)
$ is a mixingale if for finite non-negative constants $c$ and $\gamma_{t}$ and
all $s\geq0,$ $k\geq0$ it follows that $\left\Vert E\left[  \nu_{k}%
|\mathcal{M}_{0}\right]  -E\left[  \nu_{k}\right]  \right\Vert _{1}\leq
c\gamma_{k}$ and $\left\Vert \nu_{s}-E\left[  \nu_{s}|\mathcal{M}%
_{s+k}\right]  \right\Vert _{1}\leq c\gamma_{k+1}$ and $\gamma_{k}%
\rightarrow0$ as $k\rightarrow\infty$. Note in particular that the definition
in McLeish takes a slightly more complicated form in terms of the bounding
constant $c$ to account for possible non-stationarity of the process. In
addition, because here $\nu_{s}$\textbf{ }is adapted to the filtration
$\mathcal{M}_{s}$ by construction, it follows that the second condition of
McLeish (1975) is automatically satisfied since $\left\Vert \nu_{s}-E\left[
\nu_{s}|\mathcal{M}_{s+k}\right]  \right\Vert _{1}=0$.

As in Dedecker and Doukhan (2003, Lemma 2), we also assume the following:

\begin{condition}
\label{DedeckerDoukhan_Lemma2} One of the three conditions below hold:
\newline(i) $P\left(  \left\vert \nu_{s}\right\vert >x\right)  \leq\left(
c/x\right)  ^{r}$ for some $r>2$ and $\sum_{k\geq0}\left(  \gamma
_{k}/2\right)  ^{\left(  r-2\right)  /\left(  r-1\right)  }<\infty$
\newline(ii) $\left\Vert \nu_{s}\right\Vert _{r}\leq\infty$ for some $r>2$ and
$\sum_{k\geq0}k^{1/\left(  r-2\right)  }\gamma_{k}<\infty$\newline(iii)
$E\left[  \left\vert \nu_{s}\right\vert ^{2}\left(  \ln\left(  1+\left\vert
\nu_{s}\right\vert \right)  \right)  \right]  <\infty$ and $\gamma
_{k}=O\left(  a^{k}\right)  $ for some $a<1$.
\end{condition}

The mixingale type assumptions made in Condition \ref{DedeckerDoukhan_Lemma2}
represent the typical trade-off between the tail thickness of the marginal
distribution of $\nu_{t}$ and the rate of decay of the mixingale coefficients
$\gamma_{i}$. An application of Dedecker and Doukhan (2003, Lemma 2) shows
that the conditions in Condition \ref{DedeckerDoukhan_Lemma2} imply that the
condition $D\left(  2,\gamma/2,\nu\right)  $, defined in Dedecker and Doukhan
(2003, p. 71), holds. It then follows by Dedecker and Doukhan (2003,
Proposition 2) that $\operatorname*{Var}\left(  \tau^{-1/2}S_{\nu,\tau
}\right)  $ converges. In addition, by Dedecker and Doukhan (2003, Corollary
1) it follows that $\nu_{0}E\left[  S_{\nu,\tau}|\mathcal{M}_{0}\right]  $
converges in $L_{1}$. By Dedecker and Rio (2000, Theorem 1) it now follows
that $E\left[  \nu_{0}^{2}|\mathcal{I}\right]  +2E\left[  \nu_{0}S_{\nu,\tau
}|\mathcal{I}\right]  $ converges in $L_{1}$ to an $\mathcal{I}$-measurable
random variable $\eta$\textbf{, }where $\mathcal{I}$\textbf{ }denotes the
sigma-algebra of all invariant sets. If $T$ is ergodic then $\eta$ is a
constant equal to
\begin{equation}
\eta=E\left[  \nu_{0}^{2}\right]  +2\sum_{s=1}^{\infty}E\left[  \nu_{0}\nu
_{s}\right]  \label{Def_eta}%
\end{equation}
by Dedecker and Rio (2000, Remark 1). Note that (\ref{Def_eta}) is the usual
variance formula in HAC type standard errors. By Dedecker and Merlev\`{e}de
(2002, Proposition 3) it then follows that Condition \ref{DM_s2} holds with
$\eta$ defined in (\ref{Def_eta}).

\subsection{Joint Stable Convergence\label{Third}}

This section adopts an argument from Barndorff-Nielsen, Hansen, Lunde, and
Shephard (2008, BHLS) to our setting where we assume conditional independence
between $Z_{\tau}$ and $Y_{n}$. We adapt the proof of Proposition 5 (p.1524)
of BHLS to our context. To establish joint stable convergence of $Z_{\tau}$
and $Y_{n}$, we show that $\exp\left(  isZ_{\tau}+itY_{n}\right)  $ converges
weakly in $L_{1}$ for all $t$ and $s$ (see Aldous and Eagleson, 1978).

Using the results in Sections \ref{First} and \ref{Second}, and a construction
in Aldous and Eagleson (1978, p.327)\footnote{To be more specific consider the
construction of $Z.$ Aldous and Eagleson show that $Z$ can be constructed as
follows: enlarge the probability space from $\left(  \mathbb{R}^{\infty}%
\times\mathbb{R}^{\infty},\mathcal{X},P\right)  $ to $\left(  \mathbb{R}%
^{\infty}\times\mathbb{R}^{\infty}\times I,\mathcal{X\times A},P\times
\lambda\right)  $ where in $\left(  I,\mathcal{A},\lambda\right)  ,$ $I$ is
the unit interval with Lebesgue measure $\lambda$. For each $\omega$ define
$Z\left(  \omega,l\right)  $ as the inverse of the distribution function
$P\left(  Z\leq z|\mathcal{M}\right)  \left(  \omega\right)  $ and where $l$
is a uniformly distributed random variable on $\left[  0,1\right]  .$ Since in
our case, conditional on $\mathcal{M}$, $Z$ is $N\left(  0,\eta\right)  $ it
follows that $Z$ can be represented as $Z=\sqrt{\eta}\xi_{\nu}$ where
$\xi_{\nu}$ is $N\left(  0,1\right)  .$}, there are random variables $Z$ and
$Y$ on a possibly enlarged probability space such that
\[
P\left(  \left.  Z\leq z\right\vert \mathcal{A}\right)  =\int_{-\infty
}^{z/\sqrt{\eta}}g\left(  x\right)  dx\text{ a.s.}%
\]
and
\[
P\left(  \left.  Y\leq y\right\vert \mathcal{A}\right)  =\int_{-\infty
}^{y/\sqrt{\sigma^{2}\left(  \nu_{1}\right)  }}g\left(  x\right)  dx\text{
a.s.}%
\]

We establish the following main result of our paper.

\begin{theorem}
\label{Joint_CLT}$\left(  Z_{\tau},Y_{n}\right)  \rightarrow_{d}\left(
Z,Y\right)  $ ($\mathcal{A}$-stably) as $n,\tau\rightarrow\infty$ where
$Z\sim\sqrt{\eta}\xi_{\nu}$, $Y\sim\sigma\left(  \nu_{1}\right)  \xi_{y}$,
$\left(  \xi_{\nu},\xi_{y}\right)  $ is multivariate standard Gaussian with
$\operatorname*{Cov}\left(  \xi_{\nu},\xi_{y}\right)  =0$ and $\sim$ stands
for two random variables having the same distribution. In addition, $\left(
\xi_{\nu},\xi_{y}\right)  $ is independent of any $\mathcal{A}$-measurable
random variable.
\end{theorem}

\begin{proof}
See Appendix.
\end{proof}

\subsection{Extension to Short Panels}

Suppose that $T\geq2$ with $T$ fixed and finite. Let's assume the same
\textquotedblleft sufficient statistic\textquotedblright\ structure, i.e.,
\[
y_{i,t}=\Upsilon_{t}\left(  \nu_{1},\ldots,\nu_{T},u_{i}\right)  .
\]
We then have $y_{i}\equiv\left(  y_{i,1},\ldots,y_{i,T}\right)  ^{\prime}$ iid
given $\mathcal{M}_{T}$, a sigma-algebra relative to which $\left\{  \nu
_{t},t\leq T\right\}  $ is measurable. The proof of Theorem
\ref{Cross-Section-CLT} goes through without modification once $\mathcal{A}$
is replaced with $\mathcal{M}_{T}$. Let $Y_{n}\equiv\frac{1}{\sqrt{n}}%
\sum_{i=1}^{n}q\left(  y_{i}\right)  $ for some real-valued function
$q$.\textbf{ }It follows that
\[
\lim_{n\rightarrow\infty}E\left[  e^{itY_{n}}|\mathcal{M}_{T}\right]
=e^{-\frac{1}{2}t^{2}\sigma^{2}\left(  \nu_{1},...,\nu_{T}\right)  },
\]
where $\sigma^{2}\left(  \nu_{1},...,\nu_{T}\right)  =E\left[  q\left(
y_{1}\right)  ^{2}|\mathcal{M}_{T}\right]  $. As for $Z_{\tau}=\frac{1}%
{\sqrt{\tau}}\sum_{s=1}^{\tau}\nu_{s}$, note that (\ref{Timeseries_CLT})
holds\ when we replace $\mathcal{A}$ with $\mathcal{M}_{T}.$ It follows that
$Z_{\tau}\rightarrow_{d}Z$ ($\mathcal{M}_{T}$-stably) where $Z\sim\sqrt{\eta
}\xi_{\nu}$ and $\eta$ is defined as before. We can write
\[
\frac{1}{\sqrt{\tau}}\sum_{s=1}^{\tau}\nu_{s}=\frac{1}{\sqrt{\tau}}\sum
_{s=1}^{T}\nu_{s}+\frac{\sqrt{\tau-T}}{\sqrt{\tau}}\frac{1}{\sqrt{\tau-T}}%
\sum_{s=T+1}^{\tau}\nu_{s}=\frac{1}{\sqrt{\tau-T}}\sum_{s=T+1}^{\tau}\nu
_{s}+o_{p}\left(  1\right)
\]
because $\frac{1}{\sqrt{\tau}}\sum_{s=1}^{T}\nu_{s}=o_{p}\left(  1\right)  $
and $\frac{\sqrt{\tau-T}}{\sqrt{\tau}}\rightarrow1$ as $\tau\rightarrow\infty
$. Then, Dedecker and Doukhan (2003, Corollary 1), Dedecker and Rio (2000,
Theorem 1) and Proposition 3 of Dedecker and Merlev\`{e}de (2002) can be
applied to $\frac{1}{\sqrt{\tau-T}}\sum_{s=T+1}^{\tau}\nu_{s}$. Joint
convergence then follows by the same argument as in Section \ref{Third} once
$\mathcal{A}$ is understood to be $\mathcal{M}_{T}$. The joint limiting
distribution is again given by $\left(  Z_{\tau},Y_{n}\right)  \rightarrow
_{d}\left(  Z,Y\right)  $ ($\mathcal{A}$-stably).

\section{Standard Errors\label{Section_SE}}

In this section we return to the example discussed in Section \ref{Model} as
an illustration of how our joint limiting distributions can be used to obtain
limiting results for the structural and treatment parameters of interest. The
estimators for the structural parameter $\beta$ and the counterfactual causal
effect $\pi\left(  v\right)  $ are based on estimators for the in-period
treatment effect $\pi_{1}$ and the parameter governing the aggregate shock
distribution $\phi.$ We recall that $\hat{\pi}_{1}=\sum_{i=1}^{n}w_{i}%
\tilde{y}_{i,1}$\textbf{ }and $\hat{\phi}=\phi+\tau^{-1}\sum_{s=1}^{\tau}%
\nu_{s}$. Using the result in Theorem \ref{Joint_CLT} we find that the joint
distribution of $\left(  \hat{\pi}_{1},\hat{\phi}\right)  $ is
\begin{equation}
\left(  \sqrt{n}\left(  \hat{\pi}_{1}-\pi_{1}\right)  ,\sqrt{\tau}\left(
\hat{\phi}-\phi\right)  \right)  \rightarrow_{d}\left(  Y,Z\right)
(\mathcal{A}\text{-stably})\label{Example_JointLimit}%
\end{equation}
where $\mathcal{M}$ is the sigma-field generated by $\nu_{1}$ and $\left(
Y,Z\right)  \sim\left(  \sigma\left(  \nu_{1}\right)  \xi_{y},\sqrt{\eta}%
\xi_{\nu}\right)  .$ First note that
\[
\sqrt{n}\left(  \hat{\pi}_{1}-\pi_{1}\right)  =\sqrt{n}\sum_{i=1}^{n}%
w_{i}u_{i}.
\]
It is straightforward to show that $\sqrt{n}\sum_{i=1}^{n}w_{i}u_{i}%
=n^{-1/2}\sum_{i=1}^{n}y_{i,1}+o_{p}\left(  1\right)  $ with $y_{i,1}=2\left(
2d_{i,1}-1\right)  u_{i}$.\footnote{Arguments similar to the proof of Theorem
4 in Kuersteiner and Prucha (2013) justify the use of $o_{p}\left(  1\right)
$ notation in the context of stable convergence. Essentially, convergence in
probability implies joint weak convergence and thus stable convergence.}%
\textbf{ }By Theorem \ref{Cross-Section-CLT} it follows that $\sigma
^{2}\left(  \nu_{1}\right)  =4\operatorname*{Var}\left(  u_{i}\right)  $ and
$\xi_{y}\sim N\left(  0,1\right)  $. Similarly, assume that $v_{s}$ is
strictly stationary and ergodic and satisfies Condition
\ref{DedeckerDoukhan_Lemma2}. Then it follows from Lemma \ref{Lemma_DM_stable}
that $\eta=E\left[  \nu_{0}^{2}\right]  +2\sum_{s=1}^{\infty}E\left[  \nu
_{0}\nu_{s}\right]  $ and $\xi_{\nu}$ is standard Gaussian and independent of
$\xi_{y}$.

We are now in a position to derive the limiting distribution of $\hat{\beta}$
and $\hat{\pi}\left(  \nu\right)  .$ Using the expansion in (\ref{beta_expand}%
) and under the additional assumption $\sqrt{\frac{n}{\tau}}\rightarrow
\sqrt{\kappa}$ for $0<\kappa<\infty$ the limiting distribution for $\hat
{\beta}$ is given by%
\begin{equation}
\sqrt{n}\left(  \hat{\beta}-\beta\right)  \rightarrow_{d}\left(  1+\nu_{1}%
^{2}\right)  \left(  2\sigma_{u}\xi_{y}\right)  -2\pi_{1}\nu_{1}\sqrt{\kappa
}\left(  \sqrt{\eta}\xi_{\nu}\right)  \text{ }(\mathcal{A}\text{-stably}%
).\label{beta_limit}%
\end{equation}
\textbf{ }The form of the limiting distribution is mixed Gaussian with random
variance%
\[
4\left(  \left(  1+\nu_{1}^{2}\right)  ^{2}\sigma_{u}^{2}+\kappa\eta\pi
_{1}^{2}\nu_{1}^{2}\right)  ,
\]
where the randomness comes from the fact that $\nu_{1}$ is a random draw from
the stationary distribution of the process $\nu_{s}$ and that $\pi_{1}$
depends on $\nu_{1}$. By the continuous mapping theorem it follows that the
standardized statistic
\[
\Upsilon_{n}=\frac{\sqrt{n}\left(  \hat{\beta}-\beta\right)  }{2\left(
\left(  1+\nu_{1}^{2}\right)  ^{2}\sigma_{u}^{2}+\kappa\eta\pi_{1}^{2}\nu
_{1}^{2}\right)  ^{1/2}}%
\]
is asymptotically normal, i.e.,
\begin{equation}
\Upsilon_{n}\rightarrow_{d}\xi\sim N\left(  0,1\right)  .\label{beta_standard}%
\end{equation}
\textbf{ }Similarly, the limiting distribution of $\hat{\pi}\left(
\nu\right)  $ can be obtained from the expansion in (\ref{pi_expand}). Since
in our example $\hat{\pi}\left(  \nu\right)  $ is a rescaled version of
$\hat{\beta}$ the resulting limiting distribution is just a scaled version of
(\ref{beta_limit}) where the scale factor is $\left(  1+\nu^{2}\right)  $ and
where $\nu,$ as opposed to $\nu_{1}$, is a fixed constant.

Feasible inference can be based on (\ref{beta_standard}). Consider testing the
null hypothesis that $\beta=\beta_{0}$ for some fixed value $\beta_{0}.$ A
Wald test for this hypothesis is given by the t-ratio
\[
\hat{\Upsilon}_{n}=\frac{\sqrt{n}\left(  \hat{\beta}-\beta_{0}\right)
}{2\left(  \left(  1+\hat{\nu}_{1}^{2}\right)  ^{2}\hat{\sigma}_{u}^{2}%
+\kappa\hat{\eta}\hat{\pi}_{1}^{2}\hat{\nu}_{1}^{2}\right)  ^{1/2}},
\]
where $\hat{\nu}_{1}=z_{1}-\hat{\phi}$ and $\hat{\sigma}_{u}^{2}=n^{-1}%
\sum_{i=1}^{n}\hat{u}_{i}^{2}$ with $\hat{u}_{i}=\tilde{y}_{i,1}-\hat{\pi}%
_{1}d_{i,1}$. Since $\hat{\pi}_{1}-\pi_{1}=O_{p}\left(  n^{-1/2}\right)  $ and
$\hat{\phi}-\phi=O_{p}\left(  \tau^{-1/2}\right)  $ by
(\ref{Example_JointLimit}) it follows from standard arguments that $\hat{\nu
}_{1}-\nu_{1}=o_{p}\left(  1\right)  $ and $\hat{\sigma}_{u}^{2}-\sigma
_{u}^{2}=o_{p}\left(  1\right)  .$ Consistent estimators for $\eta$ can be
obtained from procedures proposed by Newey and West (1987, 1994), Andrews
(1991) or Phillips (2005) under somewhat different conditions than imposed in
Condition \ref{DedeckerDoukhan_Lemma2}. For example, Newey and West (1987) or
Andrews (1991) impose strong mixing but do not require stationarity. Since
$\alpha$-mixing processes are also mixingales by McLeish (1975), a stationary
strong mixing process can be defined in a way that it satisfies both the
conditions of Theorem \ref{Joint_CLT} and conditions in Newey and West (1987)
or Andrews (1991). Assume that $\hat{\eta}$ is a consistent estimator for
$\eta.$ It follows from the continuous mapping theorem that $\hat{\Upsilon
}_{n}-\Upsilon_{n}=o_{p}\left(  1\right)  $ and that $\hat{\Upsilon}_{n}$ has
a limiting $N\left(  0,1\right)  $ distribution under the null. As a
consequence, the null that $\beta=\beta_{0}$ is rejected against a two-sided
alternative at level $\alpha$ if $\left\vert \hat{\Upsilon}_{n}\right\vert
>c_{1-\alpha/2}$ where $c_{1-\alpha/2}$ is the $1-\alpha/2$ quantile of the
standard Gaussian distribution. Similarly, a $1-\alpha$ confidence interval
for $\hat{\beta}$ can be constructed as $\hat{\beta}\pm c_{1-\alpha/2}\left.
\left(  2\left(  \left(  1+\hat{\nu}_{1}^{2}\right)  ^{2}\hat{\sigma}_{u}%
^{2}+\kappa\hat{\eta}\hat{\pi}_{1}^{2}\hat{\nu}_{1}^{2}\right)  ^{1/2}\right)
\right/  \sqrt{n}$.

\section{Conclusions}

This paper complements results in Hahn, Kuersteiner and Mazzocco (2016) to
allow for mis-specification in the time series and cross-section models. This
is achieved by covering strictly stationary mixingales. A novel conditional
CLT for the cross-sectional data is combined with a stable CLT for stationary
mixingales to establish joint stable convergence of the combined time series
and cross-section data sets. The asymptotic variance covariance matrix for the
time series component has the familiar HAC structure. In our setting the
variance of the cross-sectional component may be a function of underlying
random shocks and thus be itself a random variable. Thus, the joint limit is
in general mixed Gaussian rather than standard normal.

Our proofs rely on strict stationary and conditional independence between
cross-sectional and time series data. Extending these results to allow for
more general heterogenous mixingales is a topic of future research. 

\newpage

\appendix{}

\section{Proofs}

\begin{proof}
[Proof of Theorem \ref{Cross-Section-CLT}]The proof closely follows the
argument in Eagleson (1975). For $P$-almost all $\omega^{\prime}\in
\mathbb{R}^{\infty}\times\mathbb{R}^{\infty}$ it follows from Lemma 1 in
Eagleson (1975) that
\[
E_{\omega^{\prime}}\left[  y_{i,1}^{2}|\mathcal{A}\right]  \left(
\omega\right)  =E\left[  y_{i,1}^{2}|\mathcal{A}\right]  \left(
\omega\right)  \text{ }Q_{\omega^{\prime}}\text{-a.s.}%
\]
Then, there is a $N_{1}\in\mathcal{X}$ such that $P\left(  N_{1}\right)  =0$
and for fixed $\omega^{\prime}\notin N_{1}$ it follows that $E_{\omega
^{\prime}}\left[  S_{n}|\mathcal{A}\right]  \left(  \omega\right)  =0$ and
$E_{\omega^{\prime}}\left[  S_{n}^{2}|\mathcal{A}\right]  \left(
\omega\right)  =nS^{2}\left(  \nu_{1}\left(  \omega\right)  \right)  $. Using
the same argument as below (6) in Eagleson (1975) we make the following
observation: since $\nu_{1}\left(  \omega\right)  $ is $\mathcal{A}%
$-measurable it follows that for fixed $\omega^{\prime}$,
\[
\nu_{1}\left(  \omega\right)  =\nu_{1}\left(  \omega^{\prime}\right)  \text{
}Q_{\omega^{\prime}}\text{-a.s.;}%
\]
in other words for a set of $\omega$ with $Q_{\omega^{\prime}}$-measure one,
the above equality holds, where the right hand side $\nu_{1}\left(
\omega^{\prime}\right)  $ is a constant. Then, for $N_{2}\in\mathcal{X}$ such
that $P\left(  N_{2}\right)  =0$ and for fixed $\omega^{\prime}\notin N_{2}$
it follows that%
\begin{equation}
\frac{1}{nS^{2}\left(  \nu_{1}\left(  \omega^{\prime}\right)  \right)  }%
\sum_{i=1}^{n}E_{\omega^{\prime}}\left[  y_{i,1}^{2}1\left\{  \left\vert
y_{i,1}\right\vert >\varepsilon\sqrt{n}\right\}  |\mathcal{A}\right]
\leq\frac{E_{\omega^{\prime}}\left[  \left\vert y_{1,1}\right\vert ^{2+\delta
}|\mathcal{A}\right]  \left(  \omega\right)  }{\sigma^{2}\left(  \nu
_{1}\left(  \omega^{\prime}\right)  \right)  \left(  \varepsilon\sqrt
{n}\right)  ^{\delta}}\rightarrow0,\ Q_{\omega^{\prime}}\text{ a.s.}%
\label{CS_CLT_D1}%
\end{equation}
In (\ref{CS_CLT_D1}) we have used the fact that, again by Lemma 1 in Eagleson
(1975),
\[
E_{\omega^{\prime}}\left[  \left\vert y_{1,1}\right\vert ^{2+\delta
}|\mathcal{A}\right]  \left(  \omega\right)  =E\left[  \left.  \left\vert
y_{1,1}\right\vert ^{2+\delta}\right\vert \mathcal{A}\right]  \left(
\omega\right)  \text{ }Q_{\omega^{\prime}}\text{-a.s.}%
\]
Since by assumption, $E\left[  \left.  \left\vert y_{1,1}\right\vert
^{2+\delta}\right\vert \mathcal{A}\right]  \left(  \omega\right)  \leq
K<\infty$ the convergence in (\ref{CS_CLT_D1}) is established $Q_{\omega
^{\prime}}$ almost surely. Then the conditions for the Lindeberg CLT hold at
least for $\omega^{\prime}$ outside a set of $P$-measure zero. Thus, by the
Lindeberg CLT, for $P$-almost all $\omega^{\prime}$ it follows that
\[
\lim_{n\rightarrow\infty}E_{\omega^{\prime}}\left[  e^{itS_{n}}|\mathcal{A}%
\right]  =\lim_{n\rightarrow\infty}E\left[  e^{itS_{n}}|\mathcal{A}\right]
=e^{-\frac{1}{2}t^{2}\sigma^{2}\left(  \nu_{1}\left(  \omega^{\prime}\right)
\right)  }.
\]

\end{proof}

\begin{proof}
[Proof of Lemma \ref{Lemma_DM_stable}]To see that Dedecker and Merlev\`{e}de
(2002, Theorem 1), hereafter DM, implies $\mathcal{A}$-stable convergence, let
$\vartheta$ be any bounded\textbf{, }$\mathcal{A}$-measurable random variable
such that $\left\vert \vartheta\right\vert <K$ for some bounded non-random
constant $K$. As in DM, let $\mathcal{H}$ be the space of continuous functions
$\varphi:\mathbb{R\rightarrow R}$ such that $\left\vert \left(  1+x^{2}%
\right)  ^{-1}\varphi\left(  x\right)  \right\vert $ is bounded. Note that
$\mathcal{H}$ contains the space of bounded continuous functions.\textbf{
}Then, $\tau^{-1/2}S_{\nu,\tau}$ converges $\mathcal{A}$ stably if for any
$\varphi\in\mathcal{H}$
\[
\lim_{n}E\left[  \left(  \varphi\left(  \tau^{-1/2}S_{\nu,\tau}\right)
-\int\varphi\left(  x\sqrt{\eta}\right)  g\left(  x\right)  dx\right)
\vartheta\right]  =0.
\]
From $\mathcal{A}\subseteq\mathcal{M}_{k}$ and the law of iterated
expectations one obtains
\begin{align*}
&  E\left[  \left(  \varphi\left(  \tau^{-1/2}S_{\nu,\tau}\right)
-\int\varphi\left(  x\sqrt{\eta}\right)  g\left(  x\right)  dx\right)
\vartheta\right]  \\
&  =E\left[  \left(  E\left[  \left.  \varphi\left(  \tau^{-1/2}S_{\nu,\tau
}\right)  -\int\varphi\left(  x\sqrt{\eta}\right)  g\left(  x\right)
dx\right\vert \mathcal{M}_{k}\right]  \right)  \vartheta\right]  \\
&  \leq KE\left[  \left\vert E\left[  \left.  \varphi\left(  \tau^{-1/2}%
S_{\nu,\tau}\right)  -\int\varphi\left(  x\sqrt{\eta}\right)  g\left(
x\right)  dx\right\vert \mathcal{M}_{k}\right]  \right\vert \right]  \\
&  =K\left\Vert E\left[  \left.  \varphi\left(  \tau^{-1/2}S_{\nu,\tau
}\right)  -\int\varphi\left(  x\sqrt{\eta}\right)  g\left(  x\right)
dx\right\vert \mathcal{M}_{k}\right]  \right\Vert _{1}.
\end{align*}
Taking the $\lim\sup$ on both sides of the inequality then establishes that
Dedecker and Merlev\`{e}de (2002, Theorem 1) implies stable convergence.
\end{proof}

\begin{proof}
[Proof of Theorem \ref{Joint_CLT}]For any bounded $\mathcal{A}$-measurable
random variable $\vartheta$ consider
\begin{align*}
E\left[  \exp\left(  isZ_{\tau}+itY_{n}\right)  \vartheta\right]   &
=E\left[  E\left[  \exp\left(  itY_{n}+isZ_{\tau}\right)  |\mathcal{A}\right]
\vartheta\right] \\
&  =E\left[  E\left[  \exp\left(  itY_{n}\right)  |\mathcal{A}\right]
E\left[  \exp\left(  isZ_{\tau}\right)  |\mathcal{A}\right]  \vartheta\right]
\\
&  =E\left[  \left(  E\left[  \exp\left(  itY_{n}\right)  |\mathcal{A}\right]
-E\left[  \exp\left(  itY\right)  |\mathcal{A}\right]  \right)  E\left[
\exp\left(  isZ_{\tau}\right)  |\mathcal{A}\right]  \vartheta\right] \\
&  +E\left[  E\left[  \exp\left(  itY\right)  |\mathcal{A}\right]  \left(
E\left[  \exp\left(  isZ_{\tau}\right)  |\mathcal{A}\right]  -E\left[
\exp\left(  isZ\right)  |\mathcal{A}\right]  \right)  \vartheta\right] \\
&  +E\left[  E\left[  \exp\left(  itY\right)  |\mathcal{A}\right]  E\left[
\exp\left(  isZ\right)  |\mathcal{A}\right]  \vartheta\right]  ,
\end{align*}
where the second equality follows from the fact that conditional on
$\mathcal{A},$ $Y_{n}$ and $Z_{\tau}$ are independent. Then, using that
$\left\vert E\left[  \exp\left(  isZ_{\tau}\right)  |\mathcal{A}\right]
\vartheta\right\vert \leq K<\infty$ a.s. it follows that
\begin{align}
&  \left\vert E\left[  \left(  E\left[  \left.  \exp\left(  itY_{n}\right)
\right\vert \mathcal{A}\right]  -E\left[  \left.  \exp\left(  itY\right)
\right\vert \mathcal{A}\right]  \right)  E\left[  \left.  \exp\left(
isZ_{\tau}\right)  \right\vert \mathcal{A}\right]  \vartheta\right]
\right\vert \nonumber\\
&  \leq KE\left[  \left\vert E\left[  \left.  \exp\left(  itY_{n}\right)
\right\vert \mathcal{A}\right]  -E\left[  \left.  \exp\left(  itY\right)
\right\vert \mathcal{A}\right]  \right\vert \right] \nonumber\\
&  =K\left\Vert E\left[  \left.  \exp\left(  itY_{n}\right)  -\exp\left(
itY\right)  \right\vert \mathcal{A}\right]  \right\Vert _{1}\rightarrow0\text{
as }n\rightarrow\infty\label{Zn_conv}%
\end{align}
by the result in Section \ref{First}, noting that $\exp\left(  itY_{n}\right)
$ and $\exp\left(  itY\right)  $ are bounded. Next note that $E\left[
\exp\left(  itY\right)  |\mathcal{A}\right]  =e^{-\frac{1}{2}t^{2}\sigma
^{2}\left(  \nu_{1}\right)  }$ a.s. is bounded. Let $\tilde{\vartheta
}=\vartheta e^{-\frac{1}{2}t^{2}\sigma^{2}\left(  \nu_{1}\right)  }$ such that
$\tilde{\vartheta}$ is $\mathcal{A}$-measurable and bounded. By stable
convergence of $Z_{\tau}$, established in Section \ref{Second}, it then
follows that
\begin{align}
&  E\left[  E\left[  \left.  \exp\left(  itY\right)  \right\vert
\mathcal{A}\right]  \left(  E\left[  \left.  \exp\left(  isZ_{\tau}\right)
\right\vert \mathcal{A}\right]  -E\left[  \left.  \exp\left(  isZ\right)
\right\vert \mathcal{A}\right]  \right)  \vartheta\right] \nonumber\\
&  =E\left[  \left(  E\left[  \exp\left(  isZ_{\tau}\right)  |\mathcal{A}%
\right]  -E\left[  \exp\left(  isZ\right)  |\mathcal{A}\right]  \right)
\tilde{\vartheta}\right]  \rightarrow0~\text{as }\tau\rightarrow\infty.
\label{Yn_conv}%
\end{align}
Together (\ref{Zn_conv}) and (\ref{Yn_conv}) imply that for $n,\tau
\rightarrow\infty,$
\begin{align*}
E\left[  \exp\left(  isZ_{\tau}+itY_{n}\right)  \vartheta\right]   &
\rightarrow E\left[  E\left[  \left.  \exp\left(  itY\right)  \right\vert
\mathcal{A}\right]  E\left[  \left.  \exp\left(  isZ\right)  \right\vert
\mathcal{A}\right]  \vartheta\right] \\
&  =E\left[  \exp\left(  itY\right)  \exp\left(  isZ\right)  \vartheta\right]
.
\end{align*}
By the Cramer-Wold theorem this implies that $\left(  Z_{\tau},Y_{n}\right)
\rightarrow\left(  Z,Y\right)  $ $\mathcal{A}$-stably as $n,\tau
\rightarrow\infty$. To identify the limits note that $E\left[  \exp\left(
itY\right)  |\mathcal{A}\right]  =E\left[  \exp\left(  it\sigma\left(  \nu
_{1}\right)  \xi_{y}\right)  |\mathcal{A}\right]  =e^{-\frac{1}{2}t^{2}%
\sigma^{2}\left(  \nu_{1}\right)  }$ where $\xi_{y}$ is $N\left(  0,1\right)
$ and therefore independent of any $\mathcal{A}$-measurable random variable.
Similarly, $E\left[  \left.  \exp\left(  isZ_{\tau}\right)  \right\vert
\mathcal{A}\right]  =E\left[  \left.  \exp\left(  is\sqrt{\eta}\xi_{\nu
}\right)  \right\vert \mathcal{A}\right]  $ where $\xi_{\nu}$ is $N\left(
0,1\right)  $ and independent of any $\mathcal{A}$-measurable random variable.
Note that the last statement follows from (\ref{Timeseries_CLT}) after setting
$\varphi\left(  .\right)  =\exp\left(  is.\right)  .$ The fact that $\xi_{\nu
}$ and $\xi_{y}$ are independent of each other follows from the assumption,
that conditional on $\mathcal{A}$, $Z_{\tau}$ and $Y_{n}$ are independent for
all $n$ and $\tau.$
\end{proof}

\end{document}